\documentstyle[11pt]{article}
\newtheorem{tm}{Theorem}[section]

\newtheorem{pr}[tm]{Proposition}
\newtheorem{rmk}[tm]{Remark}
\newtheorem{cor}[tm]{Corollary}

\newtheorem{??}[tm]{Question}

\font\tenmsb=msbm10  
\font\sevenmsb=msbm7
\font\fivemsb=msbm5
\newfam\msbfam
\textfont\msbfam=\tenmsb
\scriptfont\msbfam=\sevenmsb
\scriptscriptfont\msbfam=\fivemsb
\def\Bbb#1{{\fam\msbfam #1}}

\font\teneufm=eufm10
\font\seveneufm=eufm7
\font\fiveeufm=eufm5
\newfam\eufmfam
\textfont\eufmfam=\teneufm
\scriptfont\eufmfam=\seveneufm
\scriptscriptfont\eufmfam=\fiveeufm

\def\lorw{\longrightarrow}
\newcommand\n{\noindent}

\newcommand\ci{\cite}

\newcommand\rat{{\Bbb Q}}
\newcommand\comp{{\Bbb C}}

\newcommand\blacksquare{{\hspace*{\fill} $\Box$}}

\title{The Gysin map is compatible with mixed Hodge Structures}

\author{
Mark Andrea A.  de Cataldo\thanks{
Partially supported by N.S.F. Grant DMS 0202321 and NSA 
Grant MDA904-02-1-0100.}\, 
and Luca Migliorini\thanks{ Partially supported
by MIUR project Propriet\`a Geometriche delle Variet\`a 
Reali e Complesse and by GNSAGA.}
}

\date{}

\begin{document}\maketitle




\section {Introduction}
The language of Bivariant Theory was developed in \ci{fuma} and it 
turns out to be extremely useful 
in Intersection Theory 
and in Riemann-Roch-type questions (cfr. \ci{fulton}). 
The Chow-theoretic version 
associates a graded group 
$A^*(X \stackrel{f}{\to} Y)$ with a morphism 
$f:X \to Y$ of algebraic schemes. 

\n
There are a product structure, a proper push-forward and a pull-back.

\n
In the case
of the structural map $X \to Spec \,k$, we find  the Chow groups 
$A^*(X \stackrel{}{\to} Spec \, k) \simeq A_{-*}(X)$,  and in the case
of the identity 
map the groups $A^*(X \stackrel{Id}{\to} X)$ 
are called the Chow cohomology groups. 

\n
When $X$ is nonsingular these
latter groups 
agree with the Chow groups:
$$
A^*(X \stackrel{Id}{\to} X) \simeq A_{\dim{X}-*}(X)
$$ 
and the bivariant product structure agrees with the usual 
intersection product.

\medskip

For varieties defined over the field of complex numbers, 
the topological counterpart, i.e. Topological Bivariant Theory, 
admits the characterization
$$H^i(X \stackrel{f}{\to} Y) \simeq Hom_{D^b(Y)}(Rf_! \rat_X, \rat_Y[i]),$$ 
from which it follows immediately that
$H^*(X \stackrel{}{\to} \{point\}) \simeq Hom (H^*_c(X), \rat)=
H_{*}^{BM}(X)$, 
the Borel-Moore homology  groups, 
and $H^i(X \stackrel{Id}{\to} X)= Hom_{D^b(Y)}( \rat_X,
\rat_X[i])=H^i(X, \rat)$, 
the cohomology groups.

\medskip
It is quite natural to expect that the topological bivariant
groups can be endowed 
with natural mixed Hodge structures (MHS)
and that the maps arising in the context of the bivariant formalism should
be compatible with MHS.

\n
It seems that this question can be satisfactorily dealt-with
only  after having developed  the Hodge theory of maps, 
and in particular the Hodge-theoretic version of the Decomposition Theorem
(cfr. \ci{decatmig2}).

\medskip
In this note we shall be less ambitious and we shall deal with a very special 
yet useful case which 
does not seem to be available in the 
literature in the generality which is needed in some applications.
See, for example, those contained 
in \ci{decatmig}, where morphisms between 
cohomology groups associated with correspondences 
between global finite group quotients of smooth, not 
necessarily complete, varieties 
are stated to be maps of MHS.

\medskip
Consider a codimension $d$ regular  embedding $h:Y \lorw X$. There is 
a refined Gysin homomorphism 
$[h] \in H^{2d}(Y \stackrel {h}{\to} X) \simeq H^{2d}(X, X\setminus Y)$, 
cfr. \ci{fulton}, 19.2,
which gives, for any map $X' \lorw X$, the so-called 
refined Gysin maps $h^!: H_{*}^{BM}(X')\to  H_{*-2d}^{BM}(Y \times _X X')$. 

\medskip
Theorem  \ref{gyismhs} states  the compatibility of
the Gysin map $ h^!$ with 
the MHS 
involved and it is proved  using the  definition
of the Gysin map via specialization to 
the normal cone (\ci{fulton}, \ci{verdier}).

\bigskip
We work with algebraic varieties and algebraic schemes defined
over the field of  complex numbers.
We use cohomology etc. with rational coefficients.

\bigskip
The first-named author dedicates this paper to the memory
of Meeyoung Kim.

\section{Remarks on mixed Hodge structures}
\label{romhs}
Given a morphism $f:Y \to X$ of algebraic
schemes,
the relative cohomology groups $H^*(X \hbox{mod}\,Y)$ are given a 
MHS which is functorial in the map $f,$ in the sense that given
a commutative diagram of maps of algebraic schemes
$$
\begin{array}{ccccc}
Y' & \stackrel{i_1}\to & Y \\
\downarrow f' & & \downarrow f \\
X' & \stackrel{i_2}\to & X,
\end{array}
$$
the natural  morphism 
$H^*(X \hbox{mod}\,Y) \to H^*(X' \hbox{mod} \, Y') $ 
is a morphism of Mixed Hodge structures (cfr. \ci{ho3}, Exemple 8.3.8; see 
also \ci{du}, $\S2$). The relative cohomology groups are 
the cohomology groups of the simplicial scheme $C(f).$

\n
The  diagram above gives rise to  a map of simplicial schemes 
$C(f')\to C(f)$, hence the morphism of MHS.

\medskip
Given an algebraic scheme $U,$ one has, for every
open immersion $U \lorw U'$ into a proper algebraic scheme, $H^l_c(U) \simeq
H^l(U', U'\setminus U)$. 

\n
It follows  that
cohomology with compact supports admits a  natural
MHS which is functorial for open immersions and proper maps.

\n
Dually, since $H^{BM}_l(U) \simeq H^l_c(U)^{\vee}$, the same is true
for Borel-Moore homology.

\section{The deformation to the normal cone}
\label{tdtnc}
For more details on what follows, see \ci{fulton} and \ci{verdier}.

\medskip
Let $h:Y \lorw X$ be a closed embedding of 
algebraic schemes, and ${\cal I}_Y \subseteq {\cal O}_X$
be 
the corresponding sheaf of ideals.

\n
The {\em normal cone to $Y$ in $X$} is defined to be
the algebraic scheme 
$$
C_{Y}X \, :=\, Spec_{{\cal O}_Y} \, \bigoplus_{n\geq 0} 
{\cal I}_Y^n/{\cal I}_Y^{n+1}
$$
There are  natural maps
$ Y \lorw C_{Y}X \lorw Y,$ where 
the second one is the (affine) cone-bundle 
projection to $Y$ and
the first one 
is the closed embedding of the zero-section.

\n
If the embedding of $Y$ in $X$ is regular
of codimension $d$, then $C_Y X$ is naturally identified with
the rank $d$ normal bundle $N_{Y, X}$ of $Y$ in $X.$

\n
The exceptional divisor of the blowing-up
$Bl_Y X \lorw X$  of $X$ along $Y$
is the projectivisation $P(C_Y X).$

\medskip
Let us recall how the deformation to the normal cone  
and the canonical maps which are associated with it 
are defined.

\medskip
Let $M = Bl_{ Y \times \{0\} } X \times {\Bbb A}^1$
be the blowing-up of $X \times {\Bbb A}^1$ along
$Y \times \{0\}.$ 

\n
The exceptional divisor
can be naturally identified with $P(C_YX \oplus {\cal O}_Y).$

\n
The resulting map $\tilde{p}: M \lorw {\Bbb A}^1$ is flat.

We have that:

\n
--
the blowing-up $Bl_YX$ is naturally embedded
in $\tilde{p}^{-1}(0) \subseteq M;$

\n
--
 $\tilde{p}^{-1} (0) = P(C_Y X \oplus {\cal O}_Y )  \cup
Bl_Y X;$

\n
-- $P(C_YX \oplus {\cal O}_Y ) \cap Bl_Y X $ can be viewed
either as the exceptional divisor of the blowing-up
$Bl_Y X,$ or as the divisor at infinity of $P(C_Y X \oplus
{\cal O}_Y );$

\n
-- $C_YX$ embeds in the exceptional divisor
$P(C_YX \oplus {\cal O}_Y)$ as the complement of the divisor at infinity;

\n
-- $Y$ embeds in $C_YX \subseteq P(C_YX \oplus {\cal O}_Y)$ 
as the zero section;

\n
-- $Y \times {\Bbb A}^1$ embeds in $M,$ compatibly with the projection
to ${\Bbb A}^1;$ in particular,
 over ${\Bbb A}^1 \setminus \{0 \},$
it is just the product embedding
into $X \times {\Bbb A}^1 \setminus \{0 \};$  $Y\times \{0 \}$
is 
embedded in $p^{-1}(0)$ as the zero-section of $C_Y X \subseteq
P(C_YX \oplus O_Y).$

\medskip
Let
$$
M \setminus Bl_Y X \, =:\, M' \stackrel{p}\lorw {\Bbb A}^1
$$
be the natural flat map. The map $p$ is not proper, even if $X$ is 
complete.

\medskip
We say that the embedding $Y \subseteq X$  (i.e. the situation
at $p^{-1}(t \neq 0)$) deforms to the embedding
$Y \subseteq C_YX$ (i.e. the situation at $p^{-1}(0)$).

\n
This construction is called the {\em deformation to the normal cone}.

\medskip
There are specialization maps  (\ci{verdier}) :
$$
h^{?}: H^{BM}_l(X) \lorw H^{BM}_l(C_Y X), 
\qquad {h^{?}}^{\vee}: H^l_c(C_Y X) \lorw 
H^l_c(X) , 
$$
whose construction is recalled in the proof of \ref{gyismhs}. 

\medskip
In the case of a regular embedding
the normal cone is in fact the normal bundle and 
the flat pull-back gives an isomorphism
$H^{BM}_l(Y)   \simeq H^{BM}_{l+2d} (C_Y X)$.

\medskip
The Gysin map  $h^!$
is, by definition, the composition $H^{BM}_l(X) \lorw
H^{BM}_l(C_Y X) \lorw H^{BM}_{l-2d}(Y).$

\section{The main result}
\label{tmr}

\begin{tm}
\label{gyismhs}
Let $h:Y \lorw X$ and $C_Y X$ be as above. The natural maps
$$
h^{?}: H^{BM}_l(X) \lorw H^{BM}_l(C_Y X), \qquad 
{h^{?}}^{\vee}: H^l_c(C_Y X) \lorw 
H^l_c(X)  
$$
are maps of MHS  of type $(0,0).$

\n
If $h$ is regular of pure codimension $d,$ 
then the Gysin map and its dual
$$
h^{!} : H^{BM}_l(X) \lorw H^{BM}_{l-2d}(Y), \qquad { h^! }^{\vee}:
H^l_c(Y) \lorw H^{l+2d}_c(X) 
$$
are maps of MHS of type $(-d,-d)$ and $(d,d),$ respectively.
\end{tm}

\bigskip

\n
{\em Proof.}
We first consider a  commutative diagram of maps of algebraic schemes
$$
\begin{array}{ccccc}
U & \stackrel{j}\lorw & {\cal X} & \stackrel{i}\longleftarrow & B \\
 & p \searrow & \downarrow \pi &  \swarrow b & \\
& & S
\end{array}
$$
where $j$ is an open immersion, $B := {\cal X} \setminus U$ 
 and $\pi$ is proper.

\n
We have a distinguished triangle in
the bounded derived category   $D^b_{cc}({\cal X})$
of the category of constructible complexes of rational vector spaces
on ${\cal X}:$
$$
Rj_! j^! \rat_{\cal X} (\simeq Rj_!\rat_{U})
 \lorw \rat_{\cal X} \lorw Ri_* \rat_{B} 
\stackrel{+1}\lorw ,
$$
to which we apply $R\pi_*\simeq R\pi_!$ and get
$$
Rp_! \rat_{U} \lorw  R\pi_*\rat_{\cal X} \lorw Rb_* \rat_{B} 
\stackrel{+1}\lorw .
$$
It follows that 
$$
{\Bbb H}^l(S,Rp_! \rat_{U}) \simeq H^l({\cal X}, B), \quad \forall l.
$$

\n
Let $s \in S$ be a point. By the remarks in $\S$\ref{romhs},
the map of pairs:
$$
( \pi^{-1}(s), b^{-1}(s) ) \lorw ({\cal X}, B)
$$
induces   natural maps $H^l({\cal X}, B) \to H^l(\pi^{-1}(s), b^{-1}(s) )$
of MHS for every $l.$  Since $\pi$ is proper, 
$\pi^{-1}(s)$ is compact and the latter group
can be identified with $H^l_c(p^{-1}(s)).$ This  endows
these latter groups with  MHS.

\n
To summarize, the natural maps
$$
{\Bbb H}^l(S, Rp_! \rat_U) \to H_c^l(p^{-1}(s))
$$
are of MHS, for every $s \in S$ and for every $l.$

\medskip

\n
We apply what above  to the following situation.

\n
Let $Y' \subseteq X'$
be  an algebraic compactification of $Y \subseteq X$, 
${\cal X} := Bl_{Y'\times \{0 \} } (X' \times
{\Bbb A}^1 )$, $B:= (X' \setminus X) \times 
({\Bbb A}^1 \setminus \{0 \} ) \coprod Bl_{Y'}X' \times
\{0 \} \coprod C_{Y' \setminus Y}(X' \setminus X) ,$
$U: = {\cal X} \setminus B.$ Let $\pi:
{\cal X} \lorw {\Bbb A}^1$ be the natural map
and $p$ and $b$ the induced maps.

\n
Note that $p^{-1}( \{ 0 \} ) = C_Y X.$

\n
We have maps of MHS by choosing $s=0$ and $s =s_0 \neq 0$:
$$
H_c^{l}(p^{-1}(0) ) \stackrel{a}\longleftarrow {\Bbb H}^{l}(Rp_! \rat_U)\simeq
H^{l}€({\cal X}, B)
\lorw H_c(p^{-1}(s_0)).
$$

\n
CLAIM: $a$ is an isomorphism of MHS for every $l.$

\n
{\em Proof.} 
It is enough to prove that $a$ is an isomorphism.
Consider the fundamental system of neighborhoods of the point $\{0\} \in 
{\Bbb A}^1$ given by disks $D_r= \{ z : \, |z| <r \}.$ 
Let ${\cal X}_r:= \pi^{-1}(D_r)$, $U_r: p^{-1}(D_r),$ 
$B_r=b^{-1}(D_r)$, $p_r:= p_{|U_r}.$
Since ${\cal X}$ is a product away from $\pi^{-1}( \{ 0 \}  ),$ 
the homotopy axiom for relative cohomology ensures that
the restriction maps
${\Bbb H}^l({\Bbb A}^1, Rp_! \rat_U) \simeq
H^l({\cal X}, B) \lorw H^l({\cal X}_r, B_r)\simeq
{\Bbb H}^l(U_r, R{p_r}_! \rat_{U_r})$ are isomorphisms
for every $r>0.$

\n
The complex $Rp_! \rat_U$ is  constructible with 
respect to the stratification
$({\Bbb A}^1 \setminus \{0\}, \{ 0 \})$ of ${\Bbb A}^1$ and, 
by \ci{go-ma2}, $\S1.4$, this implies
that the restriction map 
$H^l ({\cal X}, B)\simeq {\Bbb H}^l(
{\Bbb A}^1, Rp_! \rat_U)) \stackrel{a}\lorw  H^l_c(p^{-1}(0))$
 is an isomorphism.

\medskip

\n
The first assertion of the Theorem follows.

\bigskip

\n
If the embedding $Y \subseteq X$ is regular of codimension $d$, 
then $C_Y X$ is the normal bundle
of $Y$ in $X$
and we have an isomorphism of MHS
$H^{BM}_l(Y)   \simeq H^{BM}_{l+2d} (C_Y X)$ of type $(d,d)$
and the Gysin map  $h^!$
being, by definition, the composition $H^{BM}_l(X) \lorw
H^{BM}_l(C_Y X) \lorw H^{BM}_{l-2d}(Y),$
is therefore of MHS of type $(-d,-d).$
\blacksquare

\begin{rmk}
\label{lci}
{\rm 
After simple modifications to the statement 
and to the proof,
Theorem \ref{gyismhs} holds
when $h$ is a local complete intersection  morphism.
}
\end{rmk}

The following corollary can be  used when dealing with correspondences:

\begin{cor}
\label{corrmhs}
Let $X \stackrel{p}\longleftarrow \Gamma \stackrel{q}\lorw Y$ be 
algebraic maps of algebraic varieties.
Assume that
the graph embedding $h: \Gamma \lorw \Gamma  \times X$
is regular  and that  $q$
is proper.
Then the natural map
$$
\Gamma_*: H^{BM}_{\bullet}(X) \lorw H^{BM}_{ \bullet + 
2(\dim{\Gamma} - \dim{X} )} (Y)
$$
is of MHS  of type $(\dim{\Gamma} - \dim{X},\dim{\Gamma} - \dim{X}).$
\end{cor}
{\em Proof.} We have $\Gamma_*(x) = q_* ( h^! ( [\Gamma]\times x ) )$ 
and all operations involved are compatible with MHS.
\blacksquare

\begin{rmk}
\label{rgmq}
{\rm 
The same statement above holds  when $X$ is assumed to be smooth 
or a quotient
$X=X'/G$ of a smooth variety by a finite group, 
in which case one modifies the statement
by working with $X'$ and with the
distinguished irreducible component of $\Gamma \times_X X'.$
We leave this task to the reader. 
}
\end{rmk}

\section{Examples}
    \label{ex}
The paper \ci{decatmig} contains several examples
of correspondences, stemming from
maps between surfaces, from Hilbert schemes of points on surfaces,
semi-small resolutions of singularities etc. In what follows,
for the reader's convenience, we 
offer two of the  applications of Theorem \ref{gyismhs} contained in 
\ci{decatmig}.

\subsection{Small resolutions}
\label{sr}
Let $f_{i}:X_{i} \lorw  Y,$ $i=1,\,2,$ be  {\em small} resolutions of the 
singularities of an algebraic variety $Y.$
This means that  the $X_{i}$ are nonsingular, the $f_{i}$ are proper
and birational morphisms, and the not necessarily irreducible
algebraic schemes
$X_{i}\times_{Y} X_{i} $ have exactly one component of maximal dimension 
$\dim{Y}.$ 
This component is the unique component dominating
the spaces $X_{i}$ under either projection.
See \ci{decatmig} for more on this notion and more 
references. Note that the same will be true for  
$X_{1}\times_{Y}X_{2}.$

\medskip
Let $D_{12} \subseteq X_{1}\times_{Y} X_{2}$ be the unique
irreducible component of maximal dimension $\dim{Y}.$ 

\medskip
 Theorem \ref{gyismhs} and a simple calculation using
 the calculus of correspondences
 imply the following
 
 \begin{pr}
     \label{mhsisosm}
The maps 
$$
{D_{ij}}_{*} \, : \, H^{BM}_{\bullet}(X_{i}) \lorw 
H^{BM}_{\bullet}(X_{j})
$$
are isomorphisms of MHS with inverse ${D_{ji}}_{*}.$

\n
In particular, the virtual Hodge-Deligne numbers 
of the  varieties  $X_{i},$ $i=1,\,2,$  coincide.
\end{pr}
{\em Proof.} See \ci{decatmig}, $\S$3. \blacksquare

\subsection{ Wreath products, rational double points and orbifolds}
\label{wp}
The combinatorial details of this example are rather 
lengthy. We omit them in favor of the end result.
The interested reader can see \ci{decatmig}, $\S$7.3.

Let $Y'$ be a smooth algebraic surface and $G \subseteq SL_{2}(\comp)$
be a finite group acting on $Y'$ with only isolated fixed points.

Let $Y := Y'/G$ and $f: X \lorw Y$ be its minimal resolution.

The semi-direct product  $G_{n}$ (called the Wreath
product) of $G^{n}$ with the symmetric group 
in $n$ letters $S_{n}$ acts on ${Y'}^{n}.$

The Hilbert scheme of $n-$points $X^{[n]}$ is a 
semi-small resolution 
of the singularities of ${Y'}^{n}/G_{n}.$

There is an explicit collection of  correspondences that arises in this 
situation.

There is the notion of orbifold cohomology groups
$H^{*}({Y'}^{n}/G_{n})_{orb}$
for the pair $({Y'}^{n}, G_{n}).$
These groups   carry a natural MHS.

Theorem \ref{gyismhs} allows to prove, using the 
aforementioned collection
of correspondences, the following 

\begin{pr}
    \label{orbimhs}
    There is a canonical isomorphism of MHS
    $$
    H^{*}({Y}'^{n}/G_{n})_{orb} \, \simeq \, H^{*}(X^{[n]}).
    $$
    \end{pr}







\bigskip
\noindent
Authors' addresses:

\smallskip
Mark Andrea A. de Cataldo,
Department of Mathematics,
SUNY at Stony Brook,
Stony Brook,  NY 11794, USA. \quad 
e-mail: {\em mde@math.sunysb.edu}

\smallskip
Luca Migliorini,
Dipartimento di Matematica, Universit\`a di Bologna,
Piazza di Porta S. Donato 5,
40126 Bologna,  ITALY. \quad
e-mail: {\em migliori@dm.unibo.it}

\end{document}